\documentclass[12pt]{amsart}
\usepackage{latexsym,graphicx}
\usepackage[all,cmtip]{xy}

\numberwithin{equation}{section}
\theoremstyle{plain}

\theoremstyle{remark}

\theoremstyle{definition}

\newcommand{\D}{{\mathcal D}}

\newcommand{\G}{{\mathcal G}}

\newcommand{\h}{{\mathcal H}}

\newcommand{\Q}{\mathbb Q}
\newcommand{\R}{\mathbb R}

\newcommand{\dist}{\operatorname{dist}}

\newcommand{\id}{\operatorname{id}}

\newcommand{\supp}{\operatorname{Supp}}

\def\La{\Lambda}

\def\XXint#1#2#3{{\setbox0=\hbox{$#1{#2#3}{\int}$}
      \vcenter{\hbox{$#2#3$}}\kern-.5\wd0}}

\begin{document}

\def\cal{\mathcal}




  \def\hharr#1#2{\ \smash{\mathop{\hbox to .3in{\rightarrowfill}}\limits^{\scriptstyle#1}_{\scriptstyle#2}}\ }

\def\demo#1{\noindent {\bf #1.}}
\def\proclaim#1{\noindent {\TR #1.}}
\def\tag#1{\eqno{(#1)}}
\def\subheading#1{\noindent {\bf #1} }
\def\Cal{\mathcal}
\def\bold{\mathbf}
\def\define{\def}
\define\cjoin{\sharp}
\define\csus{\Sigma \!\!\!\! /\,}
\define\Bbb{\bold}
\define\bbz{\Bbb Z}   
\define\bbq{\Bbb Q} 
\define\bbr{\Bbb R} 
\define\bbc{\Bbb C} 
\define\bbh{\Bbb H} 
\define\bbp{\Bbb P} 
\define\bbSq{{\Bbb S}{\rm q}}
\def\pc#1{\bbp_{\bbc}^{#1}}
\def\pr#1{\bbp_{\bbr}^{#1}}
\def\cyc#1#2{\Cal C^{#1}(\Bbb P^{#2})}
\def\cycd#1#2#3{\Cal C^{#1}_{#2}(\Bbb P^{#3})}

\define\plusminus{\underline{+}}
\define\bp{\bold P}
\define\bz{\bold Z}
\define\bc{\bold C}
\define\bm{\bold M}
\define\br{\bold R}
\define\bdc{\bold c}
\define\ca{\Cal A}
\define\cd{\Cal D}
\define\cf{\Cal F}
\define\co{\Cal O}
\define\cm{\Cal M}
\define\cl{\Cal L}
\define\cu{\Cal U}
\define\cz{\Cal Z}
\define\ct{\Cal T}
\define\ctd{\ct_{\cdot}}
\define\a{\alpha}
\define\b{\beta}
\define\g{\gamma}
\define\G{\Gamma}
\define\d{\delta}
\define\r{\rho}
\define\s{\sigma}
\define\z{\zeta}
\define\x{\xi}
\define\e{\epsilon}
\define\D{\Delta}
\define\k{\kappa}
\define\ztpn{{\Cal Z}_{2p}(\pn)}
\define\kofz{K(\bz }
\define\kztoq{\kofz,2)\times\kofz,4)\times\dots\times\kofz,2q)}
\define\csig{\not\Sigma}
\define\naive{na{\"\i}ve }
\define\Kahler{K{\"a}hler}
\define\equdef{\equiv}
\define\dbar{/\!\!/}\define\piinv{\pi^{-1}}
\define\pntwo{\bp^{n+2}}
\define\psit{\Psi_{tD}}
\define\dist{\rm{dist}}
\define\supp{\rm{supp}}
\define\interior{\rm{interior}}
\define\degree{\rm{degree}}

\define\regv{\rm{Reg}(V)}
\define\singv{\rm{Sing}(V)}
\define\cnone{\bc^{n+1}}
\define\pnmone{\bp^{n+m+1}}
\define\cpdx{\Cpdx}
\define\cptilone{\tilde {\cc}_{p+1}(\csus X)}
\define\cponesg{\ctil_{p+1}(\sx)^G}
\hyphenation{to-po-lo-gy}
\hyphenation{ho-mo-lo-gy}
\define\La{\Lambda}
\define\arr{\longrightarrow}
\define\cv{{\Cal V}}
\define\pt{\text{pt}}
\define\Mor{\text{Mor}}
\define\Map{\text{Map}}
\define\lhxy{L^sH^{q}(X;\, Y)}
\define\bds{{\bold s}} 
\define\Hom{\text{Hom}}
\define\rank{\text{rank}}
\define\lhx{L^sH^q(X)}
\define\NS{\text{NS}}
\define\xz{(X;\, \bbz)}
\define\xc{(X;\, \bbc)}
\define\po{\bp^1}
\define\Proj{\text{Proj}}
\define\harr{\hookrightarrow}
\define\bl{{\bold L}}
\define\f{\varphi}
\define\phib{\bar \varphi}
\define\psib{\bar \psi}
\define\id{\text{Id}}
\define\blbx{\vrule width5pt height5pt depth0pt}
\define\bspt{\Cal L _0}
\define\EM#1{K(\Bbb Z,{#1})}
\define\h#1{H_*(#1;\bbz)}
\define\ch#1{H^*(#1;\bbz)}
\define\SP#1{SP^{#1}(S^{2n})}
\define\SPP#1{SP^{#1}(\pc n)}
\define\SPC#1{SP^{#1}(\bbc^n)}
\define\sph2n{S^{2n}}
\define\eq{\cong}

\define\sq#1#2{Sq_{2{#1},2{#2}}}
\define\hsq#1#2{{\widehat {Sq}}_{2{#1},2{#2}}}

\define\bsq#1#2{{{\bold S}{\bold q}}_{2{#1},2{#2}}}
\define\hbsq#1#2{\widehat{{\bold S}{\bold q}}_{2{#1},2{#2}}}

\define\ssq#1#2{{{\Bbb S}{\rm q}}_{2{#1},2{#2}}}
\define\tsq#1#2{\ssq{#1}{#2}}
\define\hssq#1#2{\widehat{{\Bbb S}{\rm q}}_{2{#1},2{#2}}}

\define\bdy{\partial}
\define\Skew{\text{\bf Skew}}
\define\Sym{\text{\bf Sym}}
\define\DSkew{\text{\bf DSkew}}
\define\DSym{\text{\bf DSym}}
\define\td#1{\tilde D_{#1}}
\define\limd{\varinjlim_{d\to \infty}}
\define\ktq{{\Cal K}_{2q}}

\define\Q{Q}
\define\tQ{\widetilde{Q}}
\define\R{R}
\define\tR{\widetilde{R}}

\font\titlefont=cmr10 at 14 pt
\font \TR=cmr10 at 11 pt
\font \headfont=cmr10 at 13 pt
\font \fr = eufm10

\title[ALGEBRAIC CYCLES REPRESENTING COHOMOLOGY OPERATIONS]
{\titlefont ALGEBRAIC CYCLES REPRESENTING COHOMOLOGY OPERATIONS}

\date{\today}
\author{ \TR M.-L. Michelsohn}
\thanks
{ Partially supported by I.H.E.S.}

\maketitle

\begin{abstract}
In this paper we will show that certain universal homology classes
which are fundamental in topology are
algebraic.  To be specific, the products of Eilenberg-MacLane spaces 
$\ktq \equiv \EM 2 \times \EM 4 \times ... \times \EM {2q} $
have models which are limits of complex projective varieties. 
Precisely, we have $\ktq = \limd \cycd q d n$ where $\cycd q d n$ 
denotes the Chow variety of effective cycles of codimension $q$ and 
degree $d$ on $\pc n$. 
It is natural to ask which elements in the 
homology of $\ktq$ are represented by algebraic cycles in these 
approximations.  In this paper we find such representations for 
the even dimensional classes which are known as Steenrod squares  
(as well as their Pontrjagin and join products).  
These classes are dual to the cohomology classes which 
correspond to
the  basic cohomology operations also known as the Steenrod
squares. 

\end{abstract}

\vskip .4in

\centerline{\bf Table of Contents}

{
\baselineskip=15pt

\hskip 1.2in  0. \ \ Introduction

\hskip 1.2in 1. \ \ Recap of Cycle Space Basics

\hskip 1.2in 2. \ \ The Integral Homology of $\SP 2$

\hskip 1.2in 3. \ \ Cycles Representing $\bsq k n$

\hskip 1.2in 4. \ \ The Algebraic Cycles $\hssq k n$

\hskip 1.2in 5. \ \ The Main Theorem

\hskip 1.2in 6. \ \ Assembling the Classes 

\hskip 1.2in 7. \ \ Proliferation

 }

\newpage

\baselineskip=20pt  

\subheading{\S0. Introduction}

In this paper we will show that certain universal homology classes
which are fundamental in topology are
algebraic.  To be specific, the products of Eilenberg-MacLane spaces 
$\ktq \equdef \EM 2 \times \EM 4 \times ... \times \EM {2q} $
have models which are limits of complex projective varieties. 
In the simplest case we can write $\ktq  = \limd SP^d (\pc q)$
where $SP^d$ denotes the $d$-fold symmetric product 
(cf. [DT$_1$], [DT$_2$]).  Here the projections $SP(\pc n)\to SP(\pc m)$
were given explicitly in [FL].
More generally, $\ktq = \limd \cycd q d n$ where $\cycd q d n$ 
denotes the Chow variety of effective cycles of codimension $q$ and 
degree $d$ on $\pc n$ 
[L$_1$], [L$_2$].  It is natural to ask which elements in the 
homology of $\ktq$ are represented by algebraic cycles in these 
approximations.  In this paper we find such representations for 
the even dimensional classes which correspond to
the  basic cohomology operations known as the Steenrod
squares (as well as their Pontrjagin and join products).

More precisely, the Steenrod squares, $Sq^i$, are stable cohomology 
operations which go from $H^m(X;\bbz_2)$ to $H^{m+i}(X;\bbz_2)$.  They 
are stable in the sense that they commute with suspension.  
There are canonical dual 
homology operations, cf.  [H-Mc], denoted 
$\overline{Sq}_i$ defined by
$$
\langle \overline{Sq}_i \a,\b\rangle = \langle \a, Sq^i \b \rangle 
$$
where $\langle \ ,\ \rangle$ is the non-degenerate Kronecker pairing.

Precomposing with reduction  mod  2, denoted $\rho$, we have 
$\bbSq^i \equdef  Sq^i\rho$
which are also called
Steenrod squares and are stable cohomology 
operations which go from $H^m(X;\bbz)$ to $H^{m+i}(X;\bbz_2)$.   
Under the general isomorphism $H^m(X;G) \eq [X,K(G,m)]$ 
these operations correspond 
to classes in 
$[\EM m , K(\bbz_2,m+i)]$ acting by composition and therefore, in turn, 
correspond to stable cohomology classes in 
$H^{m+i}(\EM m ;\bbz_2)$.  They are denoted 
$\bbSq^i\iota_m$ 
or
 $Sq^i\rho\iota_m$ 
where $\iota_m$ 
is the fundamental class of $\EM m$.   

The basic homology classes which we show to be algebraic are 
two-torsion classes, 
$\sq k n\in H_{2k+2n}(\EM {2n};\bbz)$, which 
after reduction mod $2$ 
are sent by $\overline{Sq}_{2k}$ to the fundamental class of 
$\EM {2n}$ reduced  mod  2.
The classes $\sq k n $ are homology counterparts to the stable 
cohomology classes $\bbSq^{2k}\iota_{2n} = Sq^{2k}\rho\iota_{2n}$
in $ H^{2k+2n}(\EM {2n};\bbz_2)$.   

These classes $\sq k n \in 
H_{2k+2n}(\EM {2n};\bbz)$ are  exactly the 
classes of Cartan which correspond to the admissible sequences 
$(0,2k)$ in [C$_1$] , and also 
to the admissible sequences $(2k)$ of [C$_2$]

From these basic classes we can show that very much of the even 
homology of $\ktq$ is algebraic.
There are two natural pairings on the  varieties 
$\cycd q d n$ which approximate $\ktq$:   The addition map
$$
\cycd q d n \times \cycd q {d'} n \hharr {+}{ }\cycd q {d+d'} n
$$
and the biadditive {\it join} pairing
$$
\cycd q d n \times \cycd q {d'} {n'} \hharr\sharp
{ } \cycd q {dd'} {n+n'+1}
$$
Both are algebraic maps (morphisms of varieties).  Thus, beginning 
with the basic classes one can show that a great deal of the 
even-degree homology of $\ktq$ is algebraic.

 We will attempt to make this
paper accessible both to the algebraic topologist and to
the algebraic geometer and in so doing, make it also
accessible, for example, to the differential
geometer who may have become intrigued by related work
on the Chern map  [LM$_{1-2}$], [BLM], [LLM$_{1-4}$]. 

The results in this paper suggest the possibility of being able
to explicitly construct Steenrod operations in motivic theory or in 
morphic cohomology [FL]. However, such constructions do not follow immediately 
and  will form the subject of future research.


  \vskip .4in  
\subheading{\S1. Cycle space basics}

We begin by a quick recap of some basic facts about the
Chow spaces $\cyc q n$.  For further details see 
[LM$_1$], [L$_2$].  Let $\pi : \bbc^{n+1}- \{0\}\arr
\bbp^n$ be the projection onto complex projective
$n$-space.  Recall that an {\sl algebraic set} is a
subset $V$ of $\bbp^n$ such that $\pi^{-1}(V) \cup
\{0\}$  is the set of zeroes of a finite set of
homogeneous polynomials.  An {\sl algebraic variety} is
an {\sl irreducible} algebraic set, i.e., if $V = V_1
\cup V_2$, with $V_1, V_2$ algebraic, then $V_1 \subset
V_2$ or $V_2 \subset
V_1$.  An {\sl effective algebraic cycle} of dim
$p$ is a formal finite sum $c = \Sigma n_i V_i$ where
the $V_i$'s are algebraic varieties of dimension $p$ and the $n_i$'s are positive integers.  The
degree of the cycle  $c$ is  $\deg c = \Sigma n_i
\deg V_i$ where $\deg V_i$ is the homological degree of
$V_i$. The space of effective algebraic cycles in
$\bbp^n$ of codimension $q$ (i.e., of dimension $n-q$) and degree $d$ is itself an
algebraic variety, called a Chow variety, which we
will denote $\cycd q d n$. 
Now we can embed $\cycd q d n$ into $\cycd q {d+1}
n$ by choosing a fixed codimension $q$ plane
$\bspt$ in $\bbp^n$ and mapping a cycle $c$ in 
$\cycd q d n$
to 
$c + \bspt$ in $\cycd q {d+1} n$ .  Then by
passing to the direct limit we have
$$
\cyc q n \ \equiv\  \varinjlim \cycd q d n
$$

We now consider the {  complex join} of a
variety $V$ in $\bbp^n$ and a variety $W$ in
$\bbp^m$:  we place $\bbp^n$ and $\bbp^m$ in
general position in $\bbp^{n+m+1}$ and define
the {\bf complex join} of $V$ and $W$, denoted $V
 \,\cjoin \,   W$, to be the point-set union of all
complex lines from points in $V$ to points in
$W$.  cf.  [H ].  Now, for $\bbp^{n-1} \subset \bbp^n$ we
have $\bbp^n = \bbp^{n-1} \,\cjoin \,  p_n$ for a
point $p_n$ in $\bbp^n - \bbp^{n-1}$.           
   \vskip .4in \vfill\eject 
   

\subheading{\S2. The integral homology of
$SP^2(S^{2n})$}

We denote by $SP^k(X)$ the $k$-fold symmetric
product of $X$, namely, the space of unordered
$k$-tuples of points in $X$.  Then by choosing a fixed
base point, $x_0$, in $X$ and mapping a $k$-tuple,
$\{x_1,...,x_k\}$ to  $\{x_1,...,x_k,x_0\}$ and passing
to the direct limit we have the infinite symmetric
product of $X$, which we denote by $SP(X)$.  Note that 
$\cyc n n = SP(\bbp^n)$.

We know from Dold-Thom [DT$_2$] that there is a 
homotopy equivalence
 $$
SP(S^{2n}) \eq \EM {2n}.
\tag {2.0}
$$
We are, therefore,
interested in the homology of $\SP 2 $.  In fact, it
will be fundamental for us that the inclusion $\SP k
\subset  \SP { }$ induces an {\sl injection} 
$$
\h {\SP k  } \arr \h {\EM {2n}}
$$
for all $k$.

To understand the  classes of interest we shall
decompose $\SP 2$ as follows.  Set
$$
\D = \{\{x,x\}\in \SP 2 : x \in S^{2n}\} =
S^{2n}
$$
and
$$
\tilde \D = \{\{x,-x\}\in \SP 2 : x \in S^{2n}\} =
\pr {2n}
$$
and let $U$ be an $\e$-tubular neighborhood of $\D$ in
$\SP 2$ in the natural metric, for small $\e$.  Set
$$
V = \SP 2 - \D.
$$
\proclaim{LEMMA 2.1}  There are homotopy equivalences
$$
U \eq S^{2n}, \ \ \ \  V \eq \pr {2n},
\ \ \ \  U \cap V \eq \bbp TS^{2n}
$$
where $\bbp T\sph2n$ denotes the projectivization of
the tangent bundle of $S^{2n}$.

\demo{Proof}  Note that $\sph2n = \D \subset U$ is a
deformation retract; hence, $\sph2n \eq U$.  We claim
that $\pr {2n} = \tilde\D \subset V$ is also a
deformation retract.  To see this consider an element
$\{x,y\} \in V$.  Since $x \neq y$, we see that
neither $\{x,-y\}$ nor $\{-x,y\}$ is an antipodal
pair on $\sph2n$.  Hence each pair is joined by a
unique shortest geodesic arc.  We move $x$ toward $-y$
and $y$ toward $-x$ uniformly in time so that when $t
= 1/2$ they have each gone half the distance. 
Specifically, we set
$$
x_t =
\frac{(1-t)x-ty}{||(1-t)x-ty||},\ \       
  y_t =\frac{(1-t)y-tx}{||(1-t)y-tx||}
$$ 
for $0 \leq t \leq 1/2$.  Then $\r_t  : V
\arr V$ given by $\r_t(\{x,y\}) = \{x_t,y_t\}$ for
$0\leq t \leq 1/2$ is a smooth homotopy with $\r_0 =
\id$ and $$
\r_{1/2}(\{x,y\}) =
\left\{\frac{x-y}{||x-y||},\frac{y-x}{||x-y||}
\right\} \in \tilde\D 
$$ for all $\{x,y\}$.  This proves $\pr {2n} = \tilde
\D \eq V $.
 
 \medskip

For the last assertion we note that $U \cap V \eq
\bdy U = \bdy \hat U /\bbz_2$ where $\bdy \hat U$
is the boundary of a tubular neighborhood,  $\hat U$,
of the diagonal in $\sph2n \times \sph2n$.   The
normal bundle, $N$, to $\D \subset \sph2n \times
\sph2n$ is equivalent to the tangent bundle of
$\sph2n$, and the action induced by $\f (x,y) =(y,x)$
on $N$ is just the bundle map $v \mapsto -v$.  Since
$\hat U$ is equivariantly equivalent to $N$, we have
$U \cap V \eq \bdy U \eq \bbp T\sph2n$ as
claimed.  
\qed

 \medskip
Applying Mayer-Vietoris to the pair $(U,V)$ gives a
long exact sequence
$$
\arr H_k(\bdy U) \arr H_k(\sph2n )\oplus H_k(\pr
{2n})\arr H_k(\SP 2 ) \arr H_{k-1}(\bdy U) \arr
$$
from which one concludes that
$$
\bdy : H_{k+1}(\SP 2;\bbz) \eq H_k(\bdy
 U;\bbz)\ \ \text{for all $k > 2n$.}\tag {2.2}
$$

We now use the Serre spectral homology sequence for the
fibration
$$
    \pr {2n-1} \arr \bdy U \arr S^{2n}\tag {2.3}
$$
which has $E^2$-term $E^2 = H_*(S^{2n}; H_*(\pr
{2n-1}))
$
to conclude that
$$
H_k(\bdy U ;\bbz ) = \begin{cases} \bbz &\text{if $k=0$ or
$4n-1$}\\ 
\bbz_2 &\text{if $k=1,3,...,4n-3$ but
$k\not= 2n-1$}\\ 
\bbz_4 &\text{if $k= 2n-1$}\\
0&\text{otherwise}
\end{cases}
$$
$$
\ \ H_k(\SP 2 ;\bbz ) = \begin{cases} \bbz &\text{if $k=0,
2n$ or $4n$}\\ 
\bbz_2 &\text{if $k =2n+2,2n+4,...,4n-2$}\\ 
0&\text{otherwise}
\end{cases}\tag {2.4}
$$
(Note that the $E^2$-term is lacunary and there is only one 
differential to compute.  See [C$_2$].)

The classes which are of interest to us here are the
generators of the torsion classes in $\h {\SP 2}$.  Let
$$
\bsq k n \in H_{2(n+k)}(\SP 2;\bbz) \tag {2.5}
$$
denote the non-zero element for $0 < k < n$.
   
    \vskip .4in


\subheading{\S3. Cycles representing $\bsq k n$}

In this section we shall present explicit real algebraic
cycles which represent the 2-torsion classes (2.5).  For
this we write $S^{2n}$ as 
$$
S^{2n}= \bbc^n \cup \{\infty\} \tag {3.1}
$$
where the embedding $\bbc^n \subset S^{2n}$ is given
by the inverse of stereographic projection.  This gives 
embeddings
$$
\bbc^n  \times \bbc^n \subset S^{2n} \times
S^{2n} \tag{ 3.2}
$$
and
$$
SP^2(\bbc^n) \subset \SP 2 \tag {3.3}
$$
Let $(z,\z) = (z_1,...,z_n;\z_1,...\z_n)$ be coordinates
for $\bbc^n \times \bbc^n$, so that the unordered
pairs $\{z,\z\}$ parametrize $SP^2(\bbc^n)$.  For
each $k$, $0<k<n$, we define cycles $\ssq k n$ by setting
$$\begin{aligned}
\ssq k n &\equdef \text{the closure in $\SP 2$  of}\\
 \bigl\{\{z,\z\}\in  
SP^2(\bbc^n) &:z_{k+1}=\z_{k+1},z_{k+2}=\z_{k+2},...,z_n
=\z_n\bigr\}
\end{aligned} \tag {3.4}
$$
To see that this is real algebraic we write
$$
S^{2n} = \big\{(w,t) \in \bbc^n \times \bbr : |w|^2 + t^2 = 1
 \big\}
$$
and include $\bbc^n \harr S^{2n}$ by $z \arr
 \left (\frac{2z}{1+|z|^2},1-\frac{2}{1+|z|^2}\right)$ 
that is, by stereographic projection.  Now the cycle 
$\ssq k n$ is the quotient by $\bbz_2$ of 
$$
\bigl\{((w,t),(w',t)) \in S^{2n} \times S^{2n} : w_i = 
w_i', i = k+1, ..., n \bigr\}
$$
and is clearly a real algebraic subset of $\SP 2$ of dimension
$2(n+k)$.  It is regular and canonically oriented (by the
complex structure) outside the diagonal $\D$.  Thereby
$\tsq k n$ becomes an integral cycle in our space (as a
rectifiable current [Fed] or via triangulation of
the pair $\bigl(\SP 2,\tsq k n\bigr)$ ).

\proclaim{THEOREM 3.5} The homology class of $\tsq k
n$ is not zero.  Hence, $\tsq k n$ represents the
2-torsion class $\bsq k n \in H_{2(n+k)}(\SP 2 ;\bbz)$ defined in \S2. 
 
\demo{Proof}  It will suffice to prove that $\tsq k n$
is not homologous to zero as a cycle modulo $2$.  To do
this we produce a mod-$2$ cycle $T$ of dimension $2(n-k)$
which meets $\tsq k n$ transversely in exactly one
regular point. The cycle $T$ is defined as follows.
Note that the ''normal'' to $\Delta$ in $SP^2(S^{2n})$
is the cone on $\pr{2n-1}$ (since we are dividing by the 
flip along the diagonal).  We shall take $T$ to be a linear
subspace $\pr{2n-k}\subset   \pr{2n-1}$.
Specifically, let $\e > 0$ be fixed.  Set
$$
SN_{\e}=\{(z,-z)\in \bbc^n \times \bbc^n :
||z||=\e\}\approx S^{2n-1},
$$
and let
$$
\begin{aligned}
\bbp N_{\e}&= SN_{\e}/\bbz_2 \\
&=\{(z,-z)\in
\bbc^n \times \bbc^n : ||z||=\e\}/\bbz_2\approx
\pr {2n-1},
\end{aligned}
$$  
Note that $\tsq k n$ meets $\bbp N_{\e}$
transversally and that 
$$
\begin{aligned}
\tsq k n \cap \bbp N_{\e} 
 &= \bigl\{\{z,-z\} \in
\bbp N_{\e} : z_{k+1}=...=z_n=0\bigr\}\\
& \cong \pr{2k-1} 
 \end{aligned} \tag {3.6}
$$
We now choose $\pr {2(n-k)} = T \subset \bbp
N_{\e} =\pr {2n-1}$ to be a projective linear
subspace of dimension $2(n-k)$ which meets the
projective linear subspace (3.6) in one point. (Almost any
choice will do.)  Then $T$ is a mod-$2$ cycle of dimension
$2(n-k)$ which is contained in the regular locus of $\SP 2$
and meets $\tsq k n$ transversely in one point as claimed.

Taking the intersection number with $T$ mod $2$
defines a $\bbz_2$-valued cocycle on $\SP 2$ (cf
[Sc].  We have just seen
that this cocycle is non-zero on $\tsq k n$.  Hence
$r [\tsq k n ]\neq 0$ where $r  :
H_{2(n+k)}(\SP  2;\bbz ) \arr H_{2(n+k)}(\SP 2 
;\bbz_2 )$ is reduction mod $2$, and so $[\tsq k n ]\neq 0$
as claimed.    \qed  

\subheading{Note 3.7}  An alternative argument for Theorem
3.5 can be given along the following lines.  Consider the
manifold with boundary $V \equiv \SP 2 -
\Delta_{\e}$, where $\Delta_{\e}$ is the $\e$-tubular 
neighborhood of the diagonal, $\Delta$, in $SP^2(S^{2n})$.  
Proper intersection defines the 
non-degenerate Lefschetz duality pairing
$$
H_{\ell}(V, \partial V)\otimes H_{4n-\ell}(V)
\longrightarrow \bbz_2,
$$
and so $[\tsq k n ]\neq 0$ in 
$H_{\ell}(V, \partial V;\, \bbz_2)$.  We ``cone off''
$\partial(\tsq k n \cap V)$ down to $\Delta$ to give a
cycle in $\SP 2$ which is not zero.  This can be read
directly out of the Mayer-Vietoris sequence.

 We   recall from [St] that the natural inclusion 
$j : \SP 2 \hookrightarrow \SP{ } \cong K(\bbz, 2n)$ induces an
injection 
$$
j_* : \h {\SP 2} \ \hookrightarrow  \  \h { K(\bbz, 2n)}.
$$ 
By the computations of H. Cartan [C$_1$], [C$_2$] there is exactly one
2-torsion class in $H_{2(n+k)}( K(\bbz, 2n);\,\bbz)$, which we shall
denote by $\sq k n$.  Theorem 3.5 therefore gives us 
the following result.

\proclaim{THEOREM 3.8} \qquad $j_*\bsq k n =
\sq k n$.

\vskip .2in\noindent

\subheading{Note 3.9}  When pushed into
$H_{2(n+k)}( K(\bbz, 2n);\,\bbz_2)$, the class 
$\sq k n$ can be interpreted in terms of the homology
operations described in [HMc].  
Recall that if $\iota_{2n}
\in H^{2k}( K(\bbz, 2n);\,\bbz)$ 
denotes the fundamental cohomology class and $\rho $ is 
reduction mod 2, then the
Steenrod operation $Sq^{2k}$ corresponds to  a non-trivial element   
$Sq^{2k}\rho\iota_{2n}\in H^{2(n+k)}( K(\bbz, 2n);\,\bbz_2)$.

Let $\eta_{2n} \in  H_{2n}( K(\bbz, 2n);\,\bbz)$ be the fundamental
homology class and again let $\rho$ denote reduction mod 2.
Let $\overline{Sq_{2k}}$ 
denote the homology operation dual
to $Sq^{2k}$ as in [HMc].  
Recall that for a cohomology operation $\phi$, its dual  
$\overline {\phi}$ is defined by 
$$
\langle \overline {\phi} \a,\b \rangle = 
\langle   \a,\phi \b \rangle  \qquad\ \ \text{for all}\ \  \a ,\b
$$
where $\langle \ \ ,\ \ \rangle$ is the Kroneker pairing.  
Then we claim that $\overline{Sq_{2k}}\rho\sq k n = \rho\eta_{2n}$.  
That is, we claim that $\langle\rho\eta_{2n},\beta\rangle = 
\langle \overline{Sq_{2k}}\rho\sq k n,\beta\rangle$ all $\beta$.  
But the only non-trivial $\beta$ is $\rho\iota_{2n}$ and we have 
$\langle \rho\sq k n,Sq^{2k}\rho\iota_{2n}\rangle \neq  0 $.
   \vskip .4in


\subheading{\S 4. The algebraic cycles $\hssq k n$}

In this section we shall define (complex) algebraic cycles
in $SP^2(\pc n)$ which represent 2-torsion elements in
integral homology and push forward to the classes $\bsq
k n$ under the natural projection $SP^2(\pc n) \arr SP^2(\pc n/ \pc{n-1}) = SP^2(S^{2n})$.
  Let us fix homogeneous coordinates 
$[z]=[z_0,\dots,z_n]$ for $\pc n$ and bihomogeneous
coordinates $([z],[\z])$ for $\pc n \times \pc n$. 
We recall that the effective algebraic cycles of degree
(1,1) on $\pc n \times \pc n$ are in one-to-one
correspondence with projective classes of 
$(n+1)\times (n+1)$ complex matrices.  To each such
matrix $A = (\!( A_{ij})\!)$ is associated the divisor
$$
D_A \,=\, \left\{([z],[\z]) \in 
\pc n \times \pc n\ :\ z^tA\z = 0\right\}.
$$

The $\bbz_2$-action $([z],[\z]) \mapsto 
([\z],[z])$ on $\pc n \times \pc n$ induces an action on
the space $\cd_{(1,1)}$ of such divisors.  The fixed-point set
$\cd_{(1,1)}^{\bbz_2}$ has two connected components:
$$
\cd_{(1,1)}^{\bbz_2}\ =\
\text{\DSkew}\amalg \text{\DSym} $$
where
$$
\text{\DSkew} = \{ D_A : A \in \text{\Skew}\}  \qquad
\text{\DSym} = \{ D_A : A \in \text{\Sym}\}
$$
and
$$ 
\text{\Skew} = \{[A]\,:\,A^t = -A\} \qquad \text{and}
\qquad\text{\Sym} = \{[A]\,:\,A^t = A\}.
$$
Each divisor $D_A \in \cd_{(1,1)}^{\bbz_2}$ is
$\bbz_2$-invariant and determines a unique {\bf
reduced} divisor 
$$
\widetilde{D}_A = {{\frac 1 2 }}\x_* (D_A)
$$
on $\SPP 2$ where $\x: \pc n \times \pc n \longrightarrow
\SPP 2$ denotes the projection.  The cycles 
$\widetilde{D}_A$
for $A \in $ \Skew (or respectively $A \in$ \Sym) are all
rationally equivalent on $\SPP 2$.

\proclaim{PROPOSITION 4.1}  The cycles 
$$
\Q_{A,B} \ =\ \widetilde{D}_A - \widetilde{D}_B
$$
for $(A,B) \in \Skew \times \Sym$ represent 2-torsion
elements in the Chow group of rational equivalence classes of
algebraic cycles on $\SPP 2$.

\demo{Proof} Note that $2\Q_{A,B}= \x_*(D_A - D_B)$ and
that $D_A - D_B$ is rationally equivalent to zero on
$\pc n \times \pc n$.
\qed \medskip

We generalize this construction to higher codimension as follows.
Let $V \subset \bbc^{n+1}$ be a linear subspace of  dimension
$k$, set $ \Skew_V= \{A \in  \Skew\,:\, V\subset \text{Null}(A)\}$, 
and consider the subvariety $\Q (V) \subset  \pc n \times \pc n$ given by 
$$
\Q  (V)\ \equdef\ \bigcap_{A \in \Skew_V} D_A.
\tag{4.2}
$$
If $k = n-1$, then $\Q  (V) = D_{A}$ for some $A \in \Skew$ of rank 2.
However, for $k<n-1$, $\Q  (V)$ is not a complete intersection.  As above we
observe that $\Q  (V)$ is $\bbz_2$-invariant  and  define the reduced
cycle 
$$
\tQ  (V) = {{
\frac 1 2}}\x_*(\Q  (V)) \qquad\text{on}\  SP^2(\pc n).
\tag{4.3}
$$

We now introduce  an explicit family of   cycles of type $\Q  (V)$ and
examine systems of  equations which define them.  
Consider the subspace  $V_{k,n} \subset \bbc^{n+1}$ of  dimension $k$
given by
 $$
V_{k,n} \ = \ \left\{z \in \bbc^{n+1} \ :\ z_0=0\ \text{and}\ z_{k+1} = 
z_{k+2} = \dots = z_{n} =0 \right\},
\tag{4.4}
$$
and define
$$
\Q_{k,n} \ = \ \Q(V_{k,n}) \qquad\text{and}\qquad
\tQ_{k,n} \ = \ \tQ(V_{k,n}).
\tag{4.5}
$$
Let  $\bbc^n_0\times \bbc^n_0 \subset \pc n \times \pc n$ be the
affine chart defined in homogeneous coordinates by $z_0 = \z_0 = 1$.
Then 
$$
\Q_{k,n}^0  \equdef \Q_{k,n}\cap (\bbc^n_0\times \bbc^n_0)
\ =\ \{(z,\z) \in \bbc^n\times \bbc^n \ :\ z_j = \z_j\ \text{for} \ k+1
\leq j \leq n\}
$$
and 
$$
\Q_{k,n}  \ = \ \overline{\Q_{k,n}^0} \qquad\text{in} \ \ \pc n \times \pc n.
$$

Consider now the rational family of affine varieties
$\Q_{k,n}^0(t) \subset \bbc^n_0\times \bbc^n_0$ defined by the equations
$$
 \left\{\begin{aligned}(1-t)(z_{k+1}-\z_{k+1}) + t \, &= \, 0  \\
z_i - \z_i  \, &= \, 0  \qquad \text{for} \ \ k+2 \leq i \leq n 
\end{aligned}  
\right. \tag {4.6_t}
 $$
for $t \in \bbc$.  For $t\neq 1$, we let $\Q_{k,n} (t)$ be the closure of
$\Q_{k,n}^0(t)$ in $\pc n\times \pc n$.  It is straightforward to verify
that in bihomogeneous coordinates $([z_0,\dots,z_n],[\z_0,\dots,\z_n])$
on $\pc n\times \pc n$, $\Q_{k,n}^0(t)$ is defined by the system of
equations:
$$
 \left\{ \begin{aligned}
(1-t)(\z_0z_{k+1}-\z_{k+1}z_0) + t\z_0z_0 \, &= \, 0  \\
(1-t)(\z_iz_{k+1}-\z_{k+1}z_i) + t\z_0z_i \, &= \, 0  \\
(1-t)(\z_iz_{k+1}-\z_{k+1}z_i) + t\z_iz_0 \, &= \, 0  \\
\z_0z_i - \z_i z_0 \, &= \, 0  \qquad \text{for} \ \ k+2 \leq i \leq n  \\
\z_jz_i - \z_iz_j  \, &= \, 0  \qquad \text{for} \ \ k+2 \leq i,j \leq n,
\end{aligned} \right.
\eqno{(4.7_t)}
 $$
(Check that when $\z_0 = z_0 = 1$, equations (4.6$_t$) and
(4.7$_t$) are equivalent, and on the divisor $\z_0  z_0 = 0$, there are no
components of dimension $n+k$ provided $t\neq 1$.)

Let $\R_{k,n}\equiv \Q_{k,n} (1)\subset  \pc n\times \pc n$ be the
subvariety defined by the equations (4.7$_1$), i.e., with $t=1$.
Note that by the first equation $\z_0z_0 = 0$, we have
$$
\text{supp}\left(\R_{k,n}\right)\  \subset\ \left(\pc{n-1}\times \pc
n\right)\cup
 \left(\pc{n}\times \pc {n-1}\right).
\tag{4.8}
$$
The cycle $\R_{k,n}$ is also $\bbz_2$-invariant, and we can define the 
reduced  cycle
$$
\tR_{k,n} = {{\frac 1 2}}\x_*(\R_{k,n}) \qquad\text{on}\  SP^2(\pc
n). \tag{4.9}
$$

\subheading{Definition 4.10}  For $0<k<n$, let $\hssq k n$ be
the algebraic $(n+k)$-cycle on $\SPP 2$ defined by
$$
\hssq k n \ =\  \tQ_{k,n} - \tR_{k,n}.
$$

\proclaim{PROPOSITION 4.11}  The cycle  
$
\hssq k n 
$
represents a  2-torsion element  in the Chow group of rational equivalence
classes of algebraic cycles on $\SPP 2$. In particular, it
represents a  2-torsion element $\hbsq k n$ in $H_{2(n+k)}(\SPP 2;\,
\bbz)$. 

\demo{Proof}    Note that $2\hssq k n  = \x_*(\Q_{k,n} - \R_{k,n})$ and
that   the cycle  $\Q_{k,n} - \R_{k,n}
=\Q_{k,n}(0) -  \Q_{k,n}(1)$ is rationally
equivalent to zero on $\pc n \times \pc n$ (via the family (4.7$_t$)).
\qed  

\subheading{Note 4.12}  Let $\omega \in H^2(\pc n;\,\bbz)$ be the
canonical generator, and  let $\pi_i:\pc n \times \pc n \arr \pc n$ denote
projection onto the $i^{\text{th}}$ factor for $i=1,2$.   Set 
$\omega_i = \pi_i^*(\omega)$.  Then a straightforward calculation shows
that 
$$
{\Cal P}(\Q_{k,n}) \ = \ \left\{ 
\begin{aligned}   
\left(\omega_1\omega_2\right)^{\ell}\left(\omega_1+\omega_2\right)
   \qquad \ \
&\text{if}\ \  n-k=2\ell +1  \\  
\left(\omega_1\omega_2\right)^{\ell}\left(\omega_1^2+\omega_2^2\right)  
\qquad \ \
&\text{if}\ \  n-k=2\ell  +2
\end{aligned}\right.
$$
where  ${\Cal P} : H_{2(n+k)}(\pc n\times \pc n;\, \bbz) \hharr {\cong} { }
H^{2(n-k)}(\pc n\times \pc n;\, \bbz)$ is the  Poincar\'e duality
map.

    \vskip .4in


\subheading{\S5. The main theorem}

Fix homogeneous coordinates $[z_0,...,z_n]$ for $\pc n$
and let $\pc {n-1} = \{[z] \in \pc n : z_0 = 0 \}$. 
There is an affine chart $ {\Psi}:\bbc^n 
\overset\approx\arr \bbc_0^n \equiv\pc n-\pc {n-1}$ given by
$(z_1,...,z_n) \mapsto [1,z_1,...z_n]$.  Consider the real
analytic map 
$$
\pc n \overset\pi\arr \sph2n \tag {5.1}
$$
with the defining properties that
$$
\pi (\pc {n-1}) = \{\infty\} \tag {5.2}
$$
and $\pi$ is  represented by the identity in the
distinguished $\bbc^n$-charts, i.e., the diagram
$$
\xymatrix{
\bbc^n  \ar[r]^{\Psi} \ar[d]_{\id} &  \pc n \ar[d]^{\pi} \\
\bbc^n   \ar[r]_{\Psi_0}& \sph2n
}
 \tag {5.3}
$$
commutes, where $\Psi_0$ is the chart considered in \S3.

Taking the cartesian product gives a map
$$
\pc n \times \pc n  \ 
\overset {\pi\times \pi} \arr
\ \sph2n\times\sph2n
$$
with
$$
(\pi\times \pi)(\pc n\times\pc {n-1}\cup \pc {n-1}\times\pc
n) = \sph2n\times\{\infty\} \cup
\{\infty\}\times\sph2n .\tag {5.4}
$$
This descends to a map
$$
\SPP 2 \overset\Pi\arr \SP 2 .\tag {5.5}
$$

\proclaim{PROPOSITION 5.6}  
$$
\Pi_*(\hssq k n) = \ssq k n
$$

\demo{Proof}Restricting to our affine coordinates gives
an open dense subset (cf.(5.3))
$$
\xymatrix{
SP^2(\bbc^n)  \ar[r] \ar[d]_{\id} &  \SPP 2 \ar[d]^{\Pi} \\
SP^2(\bbc^n)     \ar[r] & \SP 2 
}
$$
where the horizontal arrows are inclusions and where
$SP^2(\bbc^n) = \{\{[z],[\z]\}:z_0=\z_0=1\}$. 
In this open set the divisor $\tQ_{k,n}$ is
defined by the equations
$$
z_i - \z_i \,=\, 0 \qquad\text{for}\ \ k+1\leq i \leq n.
\tag {5.7}
$$
(See (4.6) of (4.7) with $t=0$.)  Furthermore, in $\SPP 2$ the cycle 
$\tR_{k,n}$ is contained in the divisor $\z_0z_0=0$ (cf. (4.8)).
Consequently   in the open subset $\SPC 2$, the cycles $\hssq
k n$ and $\ssq k n$ are defined by the same  equations (5.7).  
The other  component $\tR_{k,n}$ of $\hssq k n$ is
contained in the complement of $\SPC 2$, which by (5.4) and (4.8)
is mapped into the ``axis'' $2n$-sphere in $\SP 2$.  Now the
real dimension of $\tR_{k,n}$ is $2(n+k)$, and since $k>0$, 
this implies that $\Pi_*(\tR_{k,n}) = 0$ as a real $2(n+k)$ -cycle. 
 Hence, $\Pi_*(\hssq k n ) = \Pi_*(\tQ_{k,n}) = \ssq k n$ as claimed.
\qed

 \medskip

We now consider the commutative diagram
$$
\xymatrix{
\SPP 2   \ar[r]^{\subset} \ar[d]_{\Pi} &  \SPP{ }    \ar[r]^{\cong\qquad\quad} \ar[d] &
 K(\bbz,2) \times \cdots \times  K(\bbz,2n)\ar[d]^{\text{pr}} \\
\SP 2    \ar[r]^{\subset\qquad\qquad} &  \SP{ } = SP(\pc n/\pc {n-1})  \ar[r]^{\qquad\qquad\cong} &  K(\bbz,2n)
}
$$
and let  $J:\SPP 2 \to  K(\bbz,2n)$ be the map given 
by composition in this diagram.  Then
 combining 3.5, 3.8,  4.11 and 5.6 gives the following main
result.

\proclaim{THEOREM 5.9}  The algebraic $(n+k)$-cycles
$\hssq k n$ in $\SPP 2$ represent non-zero 
$2$-torsion elements in integral homology which push forward 
under $J_*$ to the universal classes 
  $\sq k n \in H_{2(n+k)}(\EM {2n} ;\bbz)$.

   \vskip .4in


\subheading{\S6. Assembling the classes}

Let 
$
\bbc^2 \ \subset \ 
\bbc^3 \ \subset \
\bbc^4 \ \subset \ \dots
\bbc^{n+1} \ \subset \ \dots
$
be the flag defined by embedding $\bbc^i  \subset  \bbc^{j}$ as the first
$i$ coordinates.  This gives filtrations
$$\aligned
\pc 1 \ \subset \ 
\pc 2 \ &\subset \
\pc 3 \ \subset \ \dots \ \subset
\pc n  \\
SP^2 (\pc 1 )\ \subset \ 
SP^2 (\pc 2 )\ &\subset \
SP^2 (\pc 3 )\ \subset \ \dots \ \subset
SP^2 (\pc n ).
\endaligned\tag{6.1}
$$
We use these inclusions
$
j_m:SP^2(\bbp^m_{\bbc}) \subset \SPP 2
$
 to push our cycles forward (without changing
notation):
$$
\hssq k n = (j_m)_*(\hssq k m)  
$$
for $0 < k < m \leq n$. 
Equations (4.7$_t$) show directly that in $\pc n\times \pc n$ one has
$$
\Q_{k,m} \,=\, \Q_{k,n} \cap (\pc m\times \pc m)
\qquad\text{and}\qquad
\R_{k,m} \,=\, \R_{k,n} \cap (\pc m\times \pc m)
$$
from which it follows that
$$
\hssq k m \,=\, \hssq k n \cap SP^2(\pc m).   \tag{6.2}
$$

In [FL] Friedlander and Lawson define algebraic
\lq\lq projection maps"
$$
\r_m:SP(\pc n) \arr SP(\pc m) 
$$
for $0 < m \leq n$, with the following properties:

\medskip
\noindent
{({\sl i})}
$\r_m(\bbp_{\bbc}^l)
\subseteq SP(\bbp_{\bbc}^l)$ for all
$l < m$.

\smallskip
\noindent
{({\sl ii})} The map 
$$
 \psi_m : \pc m / \pc {m-1} \arr 
SP(\pc m / \pc {m-1})
$$
induced by $\r_m$ on the subquotient is the fundamental,
degree-one inclusion.  In particular, $ \psi_m
$ induces a homotopy equivalence
$$
( \psi_m)_* :SP( \pc m / \pc {m-1})
\overset \approx \arr  SP(\pc m / \pc
{m-1}). $$

Let $\widetilde\rho_m$ denote the composition 
$SP(\pc n) \overset {\r_m} \arr SP(\pc m) \arr SP( \pc m / \pc {m-1})$

\noindent
{({\sl iii})} 
The map
$$
\r =\prod_{m=1}^n \widetilde\r_m:  SP(\pc n)
 \ \arr\ \prod_{m=1}^n SP( \pc m / \pc {m-1}) \cong
\prod_{m=1}^n K(\bbz,2m)
$$
induced by these projections is a homotopy equivalence.

\medskip

It follows that the  classes $\hbsq k m = [\hssq k m] \in 
H_{2k+2m}(SP(\pc n);\,\bbz)$ satisfy
$$
(\widetilde \r_l)_*\left(  \hbsq k m \right)\ =
\begin{cases} \bsq k m\ \ \  \text{if $l=m$}
\\ 0\ \ \  \  \ \ \qquad \text{if $l\neq m$}\end{cases}  \tag{6.3}
$$
since $ \hbsq k m \subset SP(\pc m)$.

   \vskip .4in


\subheading{\S7. Proliferation}  The fact that the basic classes $\sq k
n$ are represented  by algebraic cycles implies that in fact much of the
homology of
$$
\limd SP^d(\pc q) \ \cong\ \prod_{j=1}^q K(\bbz, 2j)
$$ 
is algebraic.  To begin we consider  the algebraic suspension map
$
\csus_{n-q} : SP^d(\pc q) \arr \cycd q d n
$
which associates to a 0-cycle $c$ its join $\csus_{n-q}(c) = c\  \#
\ \pc{n-q-1}$ with a linear subspace complementary to $\pc q$
(cf. [L$_2$]).  This map induces a homotopy equivalence
$$
\csus_{n-q} : \limd SP^d(\pc q) \ 
\hharr {\cong}{} \ \limd \cycd q d n ,
$$
and since $\csus_{n-q}$ is an algebraic map, the corresponding
homology classes $\sq k n$ on $\limd \cycd q d n$ have algebraic
representatives (already at level $d=2$).

There are two algebraic pairings:  the addition map
$$
\cycd q d n \times \cycd q {d'} n  \hharr + {} \cycd q {d+d'} n
\tag{7.1}
$$
and the biadditive { join} pairing (cf. \S 2)
$$
\cycd q d n \times \cycd q {d'} {n'}\hharr {\sharp} {}
\cycd q {dd'} {n+n'+1} . \tag{7.2}
$$

The addition map (7.1) induces the standard H-space 
multiplication on $K(\bbz,2q)$ and therefore induces the standard
{\bf Pontrjagin product} in homology.  This map clearly takes pairs of
algebraic classes to algebraic classes.  

The  join mapping (7.2)   gives a way of {\bf composing} our classes.  We
define
$$
\sq {k_1}{n_1} * \dots * \sq {k_{\ell}}{n_{\ell}}
\ =\ \left[\hssq {k_1}{n_1} \#\left(\hssq {k_2}{n_2} \#
\left(\dots \#\hssq {k_{\ell}}{n_{\ell}}\right)\dots\right)\right].
\tag{7.3}
$$
{  All such classes are also algebraic.}

Thus {\bf the integral homology classes on $\prod K(\bbz,2j)$ generated
by the basic classes $\sq k n$ under the  Pontrjagin product and the join
composition are all algebraic}.

It is reasonable to conjecture that the set of classes 
$\sq {k_1}{n_1} * \dots * \sq {k_{\ell}}{n_{\ell}}$ correspond  to
the set of  admissible sequences  $(2k_1,2k_2, ... , 2k_{\ell})$
 appearing
in Cartan [C$_2$].  If so, then a great deal of the 2-torsion in the
canonical image $$
H_{2*}(K(\bbz, 2k);\,\bbz) \ \subset \  H_{2*}\Bigl(\prod_jK(\bbz,
2j);\,\bbz\Bigr) $$
for each $i$, will be algebraic.  We will, of course, not retrieve 
any elements with a "Bockstein" appearing.

In a sequel, [M], analogous results will be established for the 
$p$-torsion
where $p$ is an odd prime.

   \vskip .4in



\end{document}